\documentclass[12pt, reqno]{amsart}
\usepackage[margin=1.3in]{geometry}
\usepackage{url}
\usepackage[parfill]{parskip}
\usepackage{graphicx}
\usepackage{amsmath}
\usepackage{amssymb}
\usepackage{epstopdf}
\usepackage{enumerate}
\def\var{\mathop{\rm Var}}

\title{Contingency tables with uniformly bounded entries}
\author{Austin Shapiro}
\date{{\small February 2011}}
\address{{\small Department of Mathematics, University of Michigan, Ann Arbor, MI 48109-1043, USA}}
\thanks{{\small The author was supported in part by NSF grants DMS 0400617 and DMS 0856640.}}

\begin{document}
\newtheorem{Lem}{Lemma}
\newtheorem{Thm}{Theorem}
\newtheorem{Def}{Definition}

\begin{abstract}
We consider nonnegative integer matrices with specified row and column sums and upper bounds on the entries. We show that the logarithm of the number of such matrices is approximated by a concave function of the row and column sums. We give efficiently computable estimators for this function, including one suggested by a maximum-entropy random model; we show that these estimators are asymptotically exact as the dimension of the matrices goes to $\infty$. We finish by showing that, for $\kappa\geq 2$ and for sufficiently small row and column sums, the number of matrices with these row and column sums and with entries $\leq\kappa$ is greater by an exponential factor than predicted by a heuristic of independence.
\end{abstract}

\maketitle

\section{Main objects and results}
A {\sl contingency table} is defined as a nonnegative integer matrix with specified row and column sums. Specifically, given
$$R=(r_1,r_2,\ldots,r_m)\in{\mathbb Z}_{\geq 0}^m\quad{\rm and}\quad C=(c_1,c_2,\ldots,c_n)\in{\mathbb Z}_{\geq 0}^n$$
such that
$$r_1 + r_2 + \cdots + r_m = c_1 + c_2 + \cdots + c_n = N,$$
we denote by $\Pi(R,C)$ the set of all
$A=\big(a_{ij}\big)\in{\mathbb R}_{\geq 0}^{m\times n}$
such that
$$\sum_{j=1}^n a_{ij} = r_i~~(1\leq i\leq m)\quad{\rm and}\quad\sum_{i=1}^m a_{ij} = c_j~~(1\leq j\leq n).$$
Then $\Pi(R,C)$ is a convex polytope, and its integer points are known as the contingency tables with {\sl margins} $R$ and $C$; we denote the number of such tables by $T(R,C)$.

Given a matrix $K=\big(k_{ij}\big) \in({\mathbb Z}_{\geq 0}\cup\{\infty\})^{m\times n}$, define
$$\Pi_K(R,C):=\{A\in\Pi(R,C):~~a_{ij}\leq k_{ij}\quad\forall i,j\}.$$
Thus the integer points of $\Pi_K(R,C)$ are the contingency tables with margins $R$ and $C$, bounded entrywise by $K$. We denote the number of such tables by $T_K(R,C)$.

We will be particularly concerned with the case in which $K$ is constant, with all entries equal to $\kappa\in{\mathbb Z}_{>0}\cup\{\infty\}$, in which case we abuse the notation slightly by writing $\Pi_\kappa(R,C)$ in lieu of $\Pi_K(R,C)$; $T_\kappa(R,C)$ in lieu of $T_K(R,C)$; and so on. The integer points of $\Pi_\kappa(R,C)$ are contingency tables with margins $R$ and $C$ and entries restricted to $\{0,1,2,\ldots,\kappa\}$. When $\kappa=1$, these tables are known as {\sl 0--1 tables}. When $\kappa=\infty$, it is to be understood that $\{0,1,2,\ldots,\kappa\}$ signifies ${\mathbb Z}_{\geq 0}$, so that $\Pi_\infty(R,C):=\Pi(R,C)$.

Contingency tables arise in statistics, where they are used to represent the joint distribution of two categorical variables in a sample. In this context, one may wish to test how significantly an observed table $A$ deviates from a ``typical'' table with the observed margins. A plausible candidate for the ``typical'' table is the {\sl rank~1 table},
$$A^{\rm ind} := \left(\frac{r_ic_j}{N}\right)_{i,j}.$$
Deviation from this table tests a hypothesis of independence between the two categorical variables.

In a similar spirit, Good~\cite{Good1976} proposed the following heuristic for estimating $T(R,C)$: Consider the set of $m\times n$ nonnegative integer matrices with sum of entries equal to $N$; there are ${N+mn-1 \choose mn-1}$ such matrices. Equip this set with the uniform probability measure. Then the probability that a random sample from this set has row margin $R$ is
$${N+mn-1 \choose mn-1}^{-1} ~\prod_{i=1}^m {r_i+n-1\choose n-1},$$
while the probability that a random sample has column margin $C$ is
$${N+mn-1 \choose mn-1}^{-1} ~\prod_{j=1}^n {c_j+m-1\choose m-1}.$$
If these two events were independent, then the number of tables satisfying both constraints would be
$$I(R,C) := {N+mn-1 \choose mn-1}^{-1} ~\prod_{i=1}^m {r_i+n-1\choose n-1}~\prod_{j=1}^n {c_j+m-1\choose m-1}.$$

In~\cite{DE1985}, see also~\cite{DG1995}, Diaconis and Efron criticized the choice of $A^{\rm ind}$ as typical, and discussed alternatives to the independence hypothesis, including a hypothesis that $A$ is drawn uniformly at random from the set of tables with the observed margins. This creates the need for a different notion of the ``typical'' table, which, it turns out, may diverge dramatically from the rank~1 table (see~\cite{Barvinok2010}).

Using the theory of permanents and matrix scaling, Barvinok estimated $T(R,C)$~\cite{Barvinok2008} and $T_1(R,C)$~\cite{Barvinok2010b} and showed that $T(R,C)$ is approximately log-concave~\cite{Barvinok2007}. Barvinok also showed~\cite{Barvinok2009} that, in a certain asymptotic sense (see Definition~\ref{cloning}, ``Cloning''), $T(R,C)$ is almost always severely underestimated by $I(R,C)$. That is, for most values of $R$ and $C$, when $m$ and $n$ are large, the conditions
$$\sum_{j=1}^n a_{ij} = r_i~~(1\leq i\leq m)\quad{\rm and}\quad\sum_{i=1}^m a_{ij} = c_j~~(1\leq j\leq n)$$
are strongly positively correlated. Yet if we define
$$I_1(R,C) := {mn \choose N}^{-1} ~\prod_{i=1}^m {n\choose r_i}~\prod_{j=1}^n {m\choose c_j}$$
by analogy to $I(R,C)$, then it turns out that $I_1(R,C)$ typically {\it overestimates} $T_1(R,C)$; that is, for 0--1 tables, the conditions
$$\sum_{j=1}^n a_{ij} = r_i~~(1\leq i\leq m)\quad{\rm and}\quad\sum_{i=1}^m a_{ij} = c_j~~(1\leq j\leq n)$$
are strongly {\it negatively} correlated~\cite{Barvinok2010b} for most values of $R$ and $C$.

The aim of this paper is to extend these results to the case of general $\kappa$, and in particular to examine the transition between the positive and negative correlations described above. Although asymptotic positive correlation does not occur for any margins $(R,C)$ when $\kappa=1$, we will show that it does occur for {\it some} margins whenever $\kappa\geq 2$.

\subsection{How to read this paper}
The four principal theorems of the paper are stated in sections \ref{results1} and \ref{results2}. Each is proven in its own section (\ref{bmproof}, \ref{gfproof}, \ref{hproof}, \ref{attractionproof}). The intervening sections are thematic, each introducing a concept which will be used in the statements and proofs of the four theorems. These sections contain definitions, lemmas, and proofs of lemmas, as well as informal motivating remarks; they can be read linearly as a ``story,'' or merely scanned for their formal content. The reader in a hurry is advised to browse these sections for their definitions, then proceed directly to the proof sections, referring back to the thematic sections as needed (e.g., for the statement of lemmas cited in the proofs).

\subsection{Results: part 1}
\label{results1}
We begin by adapting a ``Brunn--Minkowski-type'' result of Barvinok~\cite{Barvinok2007} to prove approximate log-concavity for $T_K(R,C)$. Using the following definitions:
\begin{Def}
\label{omega}
For a vector or matrix $V$, let $|V|$ denote the sum of the entries of $V$.
For an integer $n\geq 0$, let $\omega(n):=\frac{n^n}{n!}$ (agreeing that $0^0=1$). For a vector or matrix $V$ with nonnegative integer entries, let $\Omega(V)$ denote the sum of $\omega(v)$ over all entries $v$ of $V$.
\end{Def}
we can state:

\begin{Thm}
\label{brunnmink}
Let $\alpha_1+\alpha_2+\cdots+\alpha_p=1~~(\alpha_1,\ldots,\alpha_p\geq 0)$. Let $R^1,\ldots,R^p\in{\mathbb Z}_{\geq 0}^m$ and $C^1,\ldots,C^p\in{\mathbb Z}_{\geq 0}^n$, such that $|R^1|=\cdots=|R^p|=|C^1|=\cdots=|C^p|=N$. Let $K\in{\mathbb Z}_{\geq 0}^{m\times n}$. Define
$$R:=\alpha_1R^1+\alpha_2R^2+\cdots+\alpha_pR^p \quad{\it and}\quad C:=\alpha_1C^1+\alpha_2C^2+\cdots+\alpha_pC^p.$$
Also, define vectors $\tilde C^1,\tilde C^2,\ldots,\tilde C^p,\tilde C\in{\mathbb Z}_{\geq 0}^n$ by
$$\tilde c^t_j := \left(\sum_{i=1}^m {k_{ij}}\right) - c^t_j \quad\quad (1\leq t\leq p,~~1\leq j\leq n)$$
and
$$\tilde c_j := \left(\sum_{i=1}^m {k_{ij}}\right) - c_j \quad\quad (= \alpha_1\tilde c^1_j + \cdots + \alpha_p\tilde c^p_j).$$
(The coordinates of vectors with uppercase names are indicated by lowercase letters: e.g., $\tilde C^t=(\tilde c^t_1,\ldots,\tilde c^t_n)$.)

Then
$$\frac{\omega(|K|)T_K(R,C)}{\Omega(R)\Omega(\tilde C)\Omega(K)} \quad\geq\quad \prod_{t=1}^p \left[\frac{T_K(R^t,C^t)}{\min\{\Omega(R^t)\Omega(\tilde C^t),~~\Omega(K)\}}\right]^{\alpha_t}.$$
\end{Thm}

This theorem is somewhat opaque in itself, due to the confounding factors $\Omega(R^t)$, $\Omega(\tilde C^t)$, etc. However, some analysis reveals that these factors typically grow more slowly than the numbers $T_K$. For a precise statement, we follow~\cite{Barvinok2009} and introduce the following asymptotic regime.

\begin{Def}[Cloning]
\label{cloning}
Let
$$R=(r_1,\ldots,r_m)\in{\mathbb Z}_{\geq 0}^m \quad{\it and}\quad C=(c_1,\ldots,c_n)\in{\mathbb Z}_{\geq 0}^n.$$
Then we define
$$R^{(s)}=(sr_1,\ldots,sr_m,~sr_1,\ldots,sr_m,~~\ldots,~~sr_1,\ldots,sr_m)$$
and
$$C^{(s)}=(sc_1,\ldots,sc_n,~sc_1,\ldots,sc_n,~~\ldots,~~sc_1,\ldots,sc_n),$$
where the number of repetitions is $s$ (thus $R^{(s)}\in{\mathbb Z}_{\geq 0}^{sm}$ and $C^{(s)}\in{\mathbb Z}_{\geq 0}^{sn}$). We refer to these vectors as the {\rm $s$-fold clonings} of $R$ and $C$.

If $K\in{\mathbb Z}_{\geq 0}^{m\times n}$, then we define $K^{(s)}$ as the $sm\times sn$ matrix of form
$$
\begin{pmatrix}
K & K & \cdots & K \\ K & K & \cdots & K \\ \vdots & \vdots & \ddots & \vdots \\ K & K & \cdots & K
\end{pmatrix}
$$
(with $s$ blocks in either direction). We call this the {\rm $s$-fold cloning} of $K$.
\end{Def}

Note that the clonings are defined so that, if $A$ is a contingency table with margins $R$ and $C$, then the $sm\times sn$ matrix
$$
\begin{pmatrix}
A & A & \cdots & A \\ A & A & \cdots & A \\ \vdots & \vdots & \ddots & \vdots \\ A & A & \cdots & A
\end{pmatrix}
$$
has margins $R^{(s)}$ and $C^{(s)}$.

Now we can state

\begin{Thm}
\label{genfunktheorem}
Let $R\in{\mathbb Z}_{>0}^m$, $C\in{\mathbb Z}_{>0}^n$, and $K\in{\mathbb Z}_{>0}^{m\times n}$. Assume that $T_K(R,C)>0$, that is, there is at least one contingency table with margins $R$ and $C$, bounded entrywise by $K$. Then
$$\lim_{s\rightarrow\infty} \frac{1}{s^2} \ln T_{K^{(s)}}(R^{(s)},C^{(s)}) \quad=\quad \ln\left(\inf_{x_1,\ldots,x_m,y_1,\ldots,y_n>0} \frac{G({\bf x},{\bf y})}{{\bf x}^R{\bf y}^C} \right),$$
where
$$G({\bf x},{\bf y}) := \prod_{i=1}^m \prod_{j=1}^n
~[1+x_iy_j+(x_iy_j)^2+\cdots+(x_iy_j)^{k_{ij}}]$$
(and ${\bf x}^R$, ${\bf y}^C$ denote $x_1^{r_1}x_2^{r_2}\cdots x_m^{r_m}$, $y_1^{c_1}y_2^{c_2}\cdots y_n^{c_n}$).
\end{Thm}

To turn this into a more usable estimate, we employ some concepts from probability theory.

\subsection{Maximum-entropy models}
\label{maxentintro}
How might we ``approximate'' the uniform distribution on $\Pi=\Pi_K(R,C)$? One na\"ive approach is to construct a random matrix with entries bounded by $K$ which satisfies the row and column constraints {\it on average}, that is, in expectation. Among all possible distributions for this random model, one compelling choice is that which achieves the maximum entropy. Such a distribution necessarily assigns equal mass to all {\it bona fide} integer points of $\Pi$, while also awarding some mass to impostors outside $\Pi$. Thus the entropy of the random model overestimates the entropy of the uniform distribution on $\Pi\cap{\mathbb Z}^{m\times n}$, and so provides an upper bound on $\ln T_K(R,C)$. Theorem~\ref{entropytheorem} expresses that this upper bound is actually a tolerably good estimate. The great advantage of the random model is that it is readily computable, with all coordinates being mutually independent.

The reader interested in the general applicability of this method to integer points of polytopes is directed to~\cite{BH2010},~\cite{S1}. The underlying ``maximum-entropy principle'' has a long history (see for example~\cite{Jaynes1957},~\cite{Jaynes1957b},~\cite{GS1985}), and we do not address its overall merits here; it is not a formal prerequisite for any of the results to follow or for their proofs, but it does lend them intuitive plausibility.

The following definitions {\it are} prerequisites for a statement of the results.

\begin{Def}
\label{truncatedgeom}
A random variable $X$ is {\rm truncated geometric} with support\break
$\{0,1,2,\ldots,\kappa\}$ if there are parameters $p\in(0,1]$ and $q\in[0,\infty)$, such that
$${\bf Pr}[X=t] = pq^t \quad\quad{\it for}~t=0,1,\ldots,\kappa.$$
For symmetry, we also say that $X$ is truncated geometric with parameters $p=0$ and $q=\infty$ if ${\bf Pr}[X=\kappa]=1$; however, in what follows, explicit treatment of this case will sometimes be left to the reader.

When $\kappa=\infty$, we have already indicated that $\{0,1,2,\ldots,\kappa\}$ is to be interpreted as ${\mathbb Z}_{\geq 0} = \{0,1,2,\ldots\}$. A random variable on this support is {\rm geometric} if there are parameters $p\in(0,1]$ and $q\in[0,1)$ (in this case necessarily satisfying $p+q=1$), such that
$${\bf Pr}[X=t] = pq^t \quad\quad{\it for}~t=0,1,2,\ldots.$$
To avoid unnecessary duplication of results, we regard this as a special case of the truncated geometric distribution.
\end{Def}

Given $\kappa\in{\mathbb Z}_{\geq 0}$ and $x\in[0,\kappa]$, or $\kappa=\infty$ and $x\in[0,\infty)$, there is a unique truncated geometric distribution with support $\{0,1,2,\ldots,\kappa\}$ and expected value equal to $x$.
\begin{Def}
We denote this distribution by $TG(x;\kappa)$, and its parameters (as in Definition~\ref{truncatedgeom}) by $p(x;\kappa)$ and $q(x;\kappa)$.
\end{Def}
The parameters $p=p(x;\kappa)$ and $q=q(x;\kappa)$ are given implicitly by the equations
\begin{align}
1 &= p(1+q+q^2+\cdots+q^\kappa), \label{implicit1}\\
x &= p(q+2q^2+\cdots+\kappa q^\kappa), \label{implicit2}
\end{align}
which, to the author's knowledge, cannot be neatly solved in general (but see section~\ref{maxent} for a discussion of the simplest cases, $\kappa=1$ and $\kappa=\infty$).

Recall that the entropy of a discrete random variable $X$ (here normalized to base $e$) is given by
$$H[X] \quad:=\quad -\sum_{x\in{\rm supp}(X)} {\bf Pr}[X=x] \ln {\bf Pr}[X=x].$$
Among all probability distributions with support in $\{0,1,2,\ldots,\kappa\}$ and given expectation $x$, the greatest entropy is achieved by $TG(x;\kappa)$, as is well-known (and may be readily proved by the method of Lagrange multipliers). Hence
\begin{Def}
Given $\kappa\in{\mathbb Z}_{\geq 0}$ and $x\in[0,\kappa]$, or $\kappa=\infty$ and $x\in[0,\infty)$, let $H^{\rm max}_\kappa(x)$ denote the entropy of $TG(x;\kappa)$.
\end{Def}
We regard $H^{\rm max}_\kappa$ as a function on $[0,\kappa]$. We break off a discussion of its formula and basic properties into section~\ref{maxent}, having said all that is needed about $H^{\rm max}_\kappa$ in order to state the paper's remaining main results.

\subsection{The independence heuristic for $\Pi_\kappa(R,C)$}
\label{indepheu}
Inspired by Good's estimate $I(R,C)$, we propose and consider an estimate $I_\kappa(R,C)$ for the number of contingency tables in $\Pi_\kappa(R,C)$ when $1\leq\kappa<\infty$. (Technically, $\kappa=0$ is allowable in all that follows, though the consequences are trivial.)
\begin{Def}[``$(\kappa+1)$-nomial coefficients''~\cite{Euler1801}]
Let $\kappa$ be a positive integer. For integers $n\geq 0$ and $0\leq r\leq n\kappa$, we denote by ${n\choose r}_\kappa$ the coefficient of $x^r$ in the polynomial expansion of $(1+x+x^2+\cdots+x^\kappa)^n$.

For integers $n\geq 0$, $r\geq 0$, we define ${n\choose r}_\infty$ to be the coefficient of $x^r$ in the power series expansion of $(1+x+x^2+\cdots)^n$.
\end{Def}
Note that ${n\choose r}_1={n\choose r}$ and ${n\choose r}_\infty={r+n-1\choose r}={r+n-1\choose n-1}$. For $\kappa\neq 1,\infty$, there is (to the author's knowledge) no comparably neat exact formula for ${n\choose r}_\kappa$.

Consider the uniform probability measure on the set of $m\times n$ matrices with entries in $\{0,1,2,\ldots,\kappa\}$ and sum of entries equal to $N$. The number of such matrices is ${mn\choose N}_\kappa$. The probability that a random sample from this set has row margin $R$ is
$${mn\choose N}_\kappa^{-1}~\prod_{i=1}^m {n\choose r_i}_\kappa;$$
the probability that a random sample has column margin $C$ is
$${mn\choose N}_\kappa^{-1}~\prod_{j=1}^n {m\choose c_j}_\kappa.$$
Thus if these two events were independent, then the number of tables satisfying both constraints would be given by
\begin{Def}
\label{indepest}
$$I_\kappa(R,C) := {mn\choose N}_\kappa^{-1} ~\prod_{i=1}^m {n\choose r_i}_\kappa~\prod_{j=1}^n {m\choose c_j}_\kappa.$$
\end{Def}
Just as $T_\kappa(R,C)$ specializes to $T(R,C)$ when $\kappa=\infty$, note that $I_\infty(R,C) = I(R,C)$.

\subsection{Results: part 2}
\label{results2}
We show the following ``log-asymptotic'' formulas for $T_\kappa(R,C)$ and $I_\kappa(R,C)$:
\begin{Thm}
\label{entropytheorem}
Let $R\in{\mathbb Z}_{\geq 0}^m$, $C\in{\mathbb Z}_{\geq 0}^n$, and $\kappa\in{\mathbb Z}_{>0}\cup\{\infty\}$. Assume $|R|=|C|=N$, and assume that $T_\kappa(R,C)>0$. Recall the definition of the cloned margins $R^{(s)}, C^{(s)}$ (Definition~\ref{cloning}). Then:
\begin{enumerate}[{\rm (i)}]
\item $$\lim_{s\rightarrow\infty} \frac{1}{s^2} \ln T_\kappa(R^{(s)},C^{(s)}) \quad=\quad \max_{Z\in\Pi_\kappa(R,C)} \sum_{i=1}^m \sum_{j=1}^n H^{\rm max}_\kappa(z_{ij}).$$
\item
\begin{align*}
\lim_{s\rightarrow\infty} \frac{1}{s^2} \ln ~&I_\kappa(R^{(s)},C^{(s)}) \quad= \\
&-mnH^{\rm max}_\kappa\left(\frac{N}{mn}\right) ~+~ n\sum_{i=1}^m H^{\rm max}_\kappa\left(\frac{r_i}{n}\right) ~+~ m\sum_{j=1}^n H^{\rm max}_\kappa\left(\frac{c_j}{m}\right).
\end{align*}
\end{enumerate}
\end{Thm}

We pause to note that the right-hand side of Theorem~\ref{entropytheorem}(i) is efficiently computable, as it is the maximum of a strictly concave function over a convex polytope; see Lemma~\ref{hmaxproperties}. This expression represents the entropy of the ``random model'' described at the beginning of section~\ref{maxentintro}.

Theorem~\ref{entropytheorem} equips us to prove the following positive correlation result:
\begin{Thm}
\label{attraction}
Continue the assumptions of Theorem~\ref{entropytheorem}; further, suppose $\kappa\geq 2$. Then there exists $\delta=\delta(\kappa)\in(0,1)$, such that if $(R,C)$ satisfy
$$\left(\max_{1\leq i\leq m} r_i\right) \left(\max_{1\leq j\leq n} c_j\right) \quad<\quad \delta\kappa N$$
then
$$\lim_{s\rightarrow\infty} \frac{1}{s^2} \ln T_\kappa(R^{(s)},C^{(s)}) \quad\geq\quad
\lim_{s\rightarrow\infty} \frac{1}{s^2} \ln I_\kappa(R^{(s)},C^{(s)}),$$
with strict inequality if neither $R$ nor $C$ is a constant vector (i.e., if it is not the case that $r_1=\cdots=r_m$ or $c_1=\cdots=c_n$).
\end{Thm}

We do not present a corresponding negative correlation result. Some commentary on the prospects for such a result can be found in section~\ref{repulsion}.

\section{Approximate log-concavity for $T_K(R,C)$ and its consequences}
\label{bm}
The following ``Brunn--Minkowski-type inequality'' is proven in~\cite{Barvinok2007}:
\begin{Thm}[Barvinok]
Let $\alpha_1+\alpha_2+\cdots+\alpha_p=1~~(\alpha_1,\ldots,\alpha_p\geq 0)$. Let $R^1,\ldots,R^p\in{\mathbb Z}_{\geq 0}^m$ and $C^1,\ldots,C^p\in{\mathbb Z}_{\geq 0}^n$, such that $|R^1|=\cdots=|R^p|=|C^1|=\cdots=|C^p|=N$. Let $W\in{\mathbb Z}_{\geq 0}^{m\times n}$. Define
$$R:=\alpha_1R^1+\alpha_2R^2+\cdots+\alpha_pR^p \quad{\it and}\quad C:=\alpha_1C^1+\alpha_2C^2+\cdots+\alpha_pC^p$$
and
$$T(R,C;W) := \sum_{A\in{\mathbb Z}^{m\times n}\cap \Pi(R,C)} \prod_{i=1}^m \prod_{j=1}^n w_{ij}^{a_{ij}}.$$
Then
$$\frac{\omega(N)T(R,C;W)}{\Omega(R)\Omega(C)} \quad\geq\quad \prod_{t=1}^p \left[\frac{T(R^t,C^t;W)}{\min\{\Omega(R^t),\Omega(C^t)\}}\right]^{\alpha_t}.$$
\end{Thm}
(See Definition~\ref{omega} for the meaning of the functions $\omega$, $\Omega$.)

The matrix $W$ may be thought of as attaching {\sl weights} to the positions of an $m\times n$ matrix. If $W$ is a 0--1 matrix, then $T(R,C;W)$ counts contingency tables with enforced zeroes in the positions $(i,j)$ such that $w_{ij}=0$. We easily prove Theorem~\ref{brunnmink} by recasting the margin and $K$-boundedness conditions as a regime of enforced zeroes (in a larger matrix).

\subsection{Proof of Theorem~\ref{brunnmink}}
\label{bmproof}
Assume the hypotheses of Theorem~\ref{brunnmink}. Define vectors ${\mathcal R}^1,{\mathcal R}^2,\ldots,{\mathcal R}^p,{\mathcal R}\in{\mathbb Z}_{\geq 0}^{m+n}$ and ${\mathcal C}\in{\mathbb Z}_{\geq 0}^{mn}$ by
\begin{align*}
{\mathcal R}^t &= (r^t_1, \ldots, r^t_m, ~\tilde c^t_1, \ldots, \tilde c^t_n), \\
{\mathcal R} &= (r_1, \ldots, r_m, ~\tilde c_1, \ldots, \tilde c_n), \\
{\mathcal C} &= (k_{11}, \ldots, k_{1n}, ~k_{21}, \ldots, k_{2n}, ~~\ldots, ~~k_{m1}, \ldots, k_{mn}).
\end{align*}
Observe that
$$|{\mathcal R}^1|=|{\mathcal R}^2|=\cdots=|{\mathcal R}^p|=|{\mathcal C}|=|K|$$
and that
$${\mathcal R} = \alpha_1{\mathcal R}^1 + \alpha_2{\mathcal R}^2 + \cdots + \alpha_p{\mathcal R}^p.$$

Define $W=\big(w_{\cdot,\cdot})$ as the $(m+n)\times(mn)$ matrix with
\begin{align*}
w_{i,(i-1)n+j} &=1 &{\rm for~all}~~i=1,\ldots,m~~{\rm and}~~j=1,\ldots,n, \\
w_{m+j,(i-1)n+j} &=1 &{\rm for~all}~~i=1,\ldots,m~~{\rm and}~~j=1,\ldots,n,
\end{align*}
and zeroes in all other positions.

Given a contingency table $A=\big(a_{ij}\big)\in\Pi_K(R,C)$, we may construct a table ${\mathcal A}=\big(a'_{\cdot,\cdot}\big)\in\Pi({\mathcal R},{\mathcal C})$ by assigning
\begin{align*}
a'_{i,(i-1)n+j} &= a_{ij} &{\rm for~all}~~i=1,\ldots,m~~{\rm and}~~j=1,\ldots,n, \\
a'_{m+j,(i-1)n+j} &= k_{ij}-a_{ij} &{\rm for~all}~~i=1,\ldots,m~~{\rm and}~~j=1,\ldots,n,
\end{align*}
and zeroes in all other positions. This conversion is easily reversed, and thus gives a bijection between tables $A\in\Pi_K(R,C)$ and tables ${\mathcal A}\in\Pi({\mathcal R},{\mathcal C})$ which have enforced zeroes in all zero positions of $W$. That is,
\begin{equation}
T_K(R,C) \quad=\quad T({\mathcal R},{\mathcal C};W). \label{distension1}
\end{equation}
Similarly,
\begin{equation}
T_K(R^t,C^t) \quad=\quad T({\mathcal R}^t,{\mathcal C};W) \label{distension2}
\end{equation}
for $t=1,\ldots,p$.

Substituting ${\mathcal R}^1,\ldots,{\mathcal R}^p,{\mathcal R}$ for $R^1,\ldots,R^p,R$ in the statement of Theorem 5, as well as ${\mathcal C}$ for $C^1,\ldots,C^p,C$ and $|K|$ for $N$, we obtain the conclusion
$$\frac{\omega(|K|)T({\mathcal R},{\mathcal C};W)}{\Omega({\mathcal R})\Omega({\mathcal C})} \quad\geq\quad \prod_{t=1}^p \left[\frac{T({\mathcal R}^t,{\mathcal C};W)}{\min\{\Omega({\mathcal R}^t),\Omega({\mathcal C})\}}\right]^{\alpha_t}.$$
Using equations \eqref{distension1} and \eqref{distension2}, and recalling the definitions of ${\mathcal R},{\mathcal R}^t,{\mathcal C}$, we rewrite the above result as
$$\frac{\omega(|K|)T_K(R,C)}{\Omega(R)\Omega(\tilde C)\Omega(K)} \quad\geq\quad \prod_{t=1}^p \left[\frac{T_K(R^t,C^t)}{\min\{\Omega(R^t)\Omega(\tilde C^t),~~\Omega(K)\}}\right]^{\alpha_t},$$
proving Theorem~\ref{brunnmink}. $\blacksquare$

{\it Remarks.}\quad The reader may notice that Theorem~\ref{brunnmink} can be stated in greater generality with only trivial modifications to the proof. For instance, $T_K(R,C)$ can be replaced by a weighted function $T_K(R,C;W)$, analogous to the function $T(R,C;W)$ in the statement of Theorem~\ref{brunnmink}; also, given
$$K=\alpha_1K^1 + \alpha_2K^2 + \cdots + \alpha_pK^p,$$
Theorem~\ref{brunnmink} remains true when each instance of $K$ on the right-hand side is replaced by $K^t$.

Theorem 5, here taken for granted, was originally derived from estimates for\break$T(R,C;W)$ obtained by Barvinok~\cite{Barvinok2008} via the theory of permanents and matrix scaling. Those estimates can be converted to the $K$-bounded setting in the manner illustrated above; however, we are able to bypass this step in the theory, as the more refined instrument of Theorem 5 is directly adaptable to our needs.

\subsection{An honestly concave proxy for $\ln T_K(R,C)$}
We now begin to assemble the ingredients for the proof of Theorem~\ref{genfunktheorem}. Throughout this section, assume $K\in{\mathbb Z}_{>0}^{m\times n}$. We define a function which ``smooths over'' $\ln T_K(R,C)$:

\begin{Def}
\label{f}
For $R\in{\mathbb R}_{\geq 0}^m$, $C\in{\mathbb R}_{\geq 0}^n$, and $K\in{\mathbb Z}_{>0}^{m\times n}$, let
$$f(R,C) = f_K(R,C) := \max_{\substack{\alpha_1,\ldots,\alpha_p\geq 0 \\ \alpha_1+\cdots+\alpha_p=1 \\ \alpha_1R^1+\cdots+\alpha_pR^p=R \\ \alpha_1C^1+\cdots+\alpha_pC^p=C}} \sum_{t=1}^p \alpha_t \ln T_K(R^t,C^t).$$
(To be clear, the maximum is taken over choices of $p\geq 1$, $\alpha_1,\ldots,\alpha_p$, $R^1,\ldots,R^p$, and $C^1,\ldots,C^p$ which satisfy the indicated constraints, and for which the summation on the right is defined. If the maximum is taken over an empty set, then we regard it as $-\infty$.)
\end{Def}
Note that the maximum in Definition~\ref{f} is well-defined (allowing for the $-\infty$ case), because there are finitely many pairs $(R,C)$ for which $T_K(R,C)>0$. It is redundant to allow any repetition among $R^1,\ldots,R^p$ or $C^1,\ldots,C^p$, so the summation on the right takes on finitely many values.

\begin{Lem}[Properties of $f(R,C)$] ~
\label{fproperties}

\begin{enumerate}[{\rm (i)}]
\item $f(R,C) \geq \ln T_K(R,C)$.
\item $f$ is concave.
\item The domain of $f$ (i.e., where $f>-\infty$) is a subset of $\Pi_K(R,C)$.
\end{enumerate}
\end{Lem}

The proof of this lemma is straightforward, but the notation is cumbersome, so it is deferred to an appendix.

\begin{Lem}[Quality of approximation]
\label{fquality}
Suppose $R\in{\mathbb Z}_{>0}^m$, $C\in{\mathbb Z}_{\geq 0}^n$, and $K\in{\mathbb Z}_{>0}^{m\times n}$. Define $\tilde C=(\tilde c_1,\ldots,\tilde c_n)$ by
$\tilde c_j := \left(\sum_{i=1}^m k_{ij}\right) - c_j,$
and suppose that $\tilde C\in{\mathbb Z}_{>0}^n$. Then
\begin{align*}
f_K(R,C)-\ln T_K(R,C) \quad\leq\quad &-\ln\sqrt{2\pi|K|} + \sum_{i=1}^m \ln\sqrt{2\pi r_i} + \sum_{j=1}^n \ln\sqrt{2\pi\tilde c_j} \\
&+ (m+n)\ln \left(\frac{e}{\sqrt{2\pi}}\right).
\end{align*}
\end{Lem}

{\it Proof.}\quad By Stirling's formula,
\begin{equation}
n-\ln\sqrt{2\pi n}-\ln\left(\frac{e}{\sqrt{2\pi}}\right) \quad\leq\quad \ln\omega(n) \quad\leq\quad n-\ln\sqrt{2\pi n} \label{omegabounds}
\end{equation}
for $n\geq 1$.

Choose $\alpha_1,\ldots,\alpha_p,R^1,\ldots,R^p,C^1,\ldots,C^p$ which achieve the maximum in Definition~\ref{f}. Now apply Theorem~\ref{brunnmink} and \eqref{omegabounds}:
\begin{align*}
f_K(R,&C)-\ln T_K(R,C) ~~\leq~~
\ln\left[\frac{\omega(|K|)}{\Omega(R)\Omega(\tilde C)\Omega(K)}\cdot \prod_{t=1}^p\left(\min\{\Omega(R^t)\Omega(\tilde C^t),~~\Omega(K)\}\right)^{\alpha_t}\right] \\
&\leq\quad \ln\left[\frac{\omega(|K|)}{\Omega(R)\Omega(\tilde C)\Omega(K)}\cdot \prod_{t=1}^p \Omega(K)^{\alpha_t}\right] \\
&=\quad \ln \frac{\omega(|K|)}{\Omega(R)\Omega(\tilde C)} \\
&=\quad \ln\omega(|K|) -\sum_{i=1}^m \ln\omega(r_i) -\sum_{j=1}^n \ln\omega(\tilde c_j) \\
&\leq\quad |K|-\ln\sqrt{2\pi|K|} -\sum_i \left(r_i-\ln\sqrt{2\pi r_i}-\ln\left(\frac{e}{2\pi}\right)\right) \\
&\hskip 100 pt -\sum_j \left(\tilde c_j-\ln\sqrt{2\pi\tilde c_j}-\ln\left(\frac{e}{2\pi}\right)\right) \\
&\leq\quad -\ln\sqrt{2\pi|K|} + \sum_{i=1}^m \ln\sqrt{2\pi r_i} + \sum_{j=1}^n \ln\sqrt{2\pi\tilde c_j}
+(m+n)\ln \left(\frac{e}{\sqrt{2\pi}}\right).\quad\square
\end{align*}

\subsection{Exact and approximate generating functions for tables}
Theorem~\ref{genfunktheorem} refers to the following polynomial:
$$G({\bf x},{\bf y}) := \prod_{i=1}^m \prod_{j=1}^n ~[1+x_iy_j+(x_iy_j)^2+\cdots+(x_iy_j)^{k_{ij}}]$$
(${\bf x}=(x_1,\ldots,x_m),~{\bf y}=(y_1,\ldots,y_n)$).

\begin{Lem}
$G$ is a generating function for $K$-bounded contingency tables; that is,
\begin{equation}
G({\bf x},{\bf y}) \quad=\quad \sum_R \sum_C T_K(R,C) {\bf x}^R{\bf y}^C, \label{genfunk}
\end{equation}
where the sum is taken over all possible margins $R$, $C$ (of lengths $m$ and $n$).
\end{Lem}

{\it Proof.}\quad Trivial. $\square$

In principle, we can ``compute'' $T_K(R,C)$ by expanding $G({\bf x},{\bf y})$ and extracting the coefficient of ${\bf x}^R{\bf y}^C$. This is of course not practical, but we might estimate this coefficient by
$$\inf_{x_i,y_j>0} \frac{G({\bf x},{\bf y})}{{\bf x}^R{\bf y}^C};$$

indeed, this is an {\it upper} bound on $T_K(R,C)$, as may be readily seen by dividing both sides of \eqref{genfunk} by ${\bf x}^R{\bf y}^C$. To bound $T_K(R,C)$ from the other side, we replace $G({\bf x},{\bf y})$ by an approximate version with smoother coefficients:

\begin{Def}
Let
$$\tilde G({\bf x},{\bf y}) \quad:=\quad \sum_R\sum_C e^{f(R,C)}{\bf x}^R{\bf y}^C,$$
where the sum is taken over all integer margins $(R,C)$ such that $f(R,C)>-\infty$.
\end{Def}
(See Definition~\ref{f} for the meaning of $f(R,C)$.)

We will find the following lemma useful, as it will allow us to pick out any nonzero term of $\tilde G({\bf x},{\bf y})$ as the largest:

\begin{Lem}
\label{tilt}
For any $(R_*,C_*)$ in the relative interior of the domain of $f$, there exist ${\bf x}_*, {\bf y}_*>0$ such that the function
$$\Phi(R,C) \quad:=\quad e^{f(R,C)}{\bf x}_*^R{\bf y}_*^C$$
attains its maximum at $R=R_*$, $C=C_*$.
\end{Lem}

{\it Proof.}\quad
Recall that $f$ is concave; therefore, its graph has a supporting hyperplane over $(R_*,C_*)$. Let such a hyperplane have outward-pointing normal vector $(u_1,\ldots,u_m,v_1,\ldots,v_n,1)$. Set
$${\bf x}_* = (x_1,\ldots,x_m) = (e^{-u_1},\ldots,e^{-u_m})\quad{\rm and}\quad
{\bf y}_* = (y_1,\ldots,y_n) = (e^{-v_1},\ldots,e^{-v_n}).$$
Then
$$\phi(R,C) := f(R,C) + \sum_{i=1}^m r_i\ln x_i + \sum_{j=1}^n c_j\ln y_j$$
is concave with respect to $R$ and $C$, and attains a critical point (hence its global maximum) at $(R_*,C_*)$. Therefore, so does $\Phi(R,C) = e^{\phi(R,C)}$. $\square$

\subsection{Proof of Theorem~\ref{genfunktheorem}}
\label{gfproof}
Assume the hypotheses of Theorem~\ref{genfunktheorem}. Using Lemma~\ref{tilt}, choose ${\bf x}_*,{\bf y}_*$ so that $e^{f(R,C)}{\bf x}^R{\bf y}^C$ is the largest term in the expansion of $\tilde G({\bf x},{\bf y})$, evaluated at ${\bf x}={\bf x}_*$ and ${\bf y}={\bf y}_*$. Thus
$$\frac{\tilde G({\bf x}_*,{\bf y}_*)}{{\bf x}_*^R{\bf y}_*^C} \quad\leq\quad
{\rm [\#~of~terms~of~}\tilde G{\rm ~with~nonzero~coeffs.]} \cdot e^{f(R,C)}.$$
The number of terms of $\tilde G$ is at most
$${\mathcal N} := \prod_{i=1}^m \left(1+\sum_{j=1}^n k_{ij}\right) + \prod_{j=1}^n \left(1+\sum_{i=1}^m k_{ij}\right),$$
since $T_K(R,C)>0$ implies that $R$ and $C$ do not exceed the margins of $K$.

Let the symbol $\heartsuit$ denote the quantity
$$-\ln\sqrt{2\pi|K|} + \sum_{i=1}^m \ln\sqrt{2\pi r_i} + \sum_{j=1}^n \ln\sqrt{2\pi\tilde c_j}
+ (m+n)\ln \left(\frac{e}{\sqrt{2\pi}}\right),$$
last seen in Lemma~\ref{fquality}.

We deduce the following chain of inequalities:
\begin{align}
\ln\left(\inf_{x_i,y_j>0} \frac{G({\bf x},{\bf y})}{{\bf x}^R{\bf y}^C} \right)
\quad&\geq\quad \ln T_K(R,C) \nonumber\\
&\geq\quad f(R,C)-\heartsuit \nonumber\\
&\geq\quad \ln\left(\frac{\tilde G({\bf x}_*,{\bf y}_*)}{{\bf x}_*^R{\bf y}_*^C \cdot{\mathcal N}}\right) - \heartsuit \nonumber\\
&\geq\quad \ln\left(\inf_{x_i,y_j>0} \frac{\tilde G({\bf x},{\bf y})}{{\bf x}^R{\bf y}^C \cdot{\mathcal N}} \right) - \heartsuit \nonumber\\
&\geq\quad \ln\left(\inf_{x_i,y_j>0} \frac{G({\bf x},{\bf y})}{{\bf x}^R{\bf y}^C} \right) - \ln{\mathcal N} - \heartsuit. \label{errorbound}
\end{align}

Now we consider the cloning of the margins. Let $G^{(s)}$ denote the generating function for $K^{(s)}$-bounded contingency tables. Letting
\begin{align*}
{\bf x}^{(s)} &:= ({\bf x}^1, {\bf x}^2, \ldots, {\bf x}^s) \\
&= (x^1_1, \ldots, x^1_m, ~x^2_1, \ldots, x^2_m, ~~\ldots, ~~x^s_1, \ldots, x^s_m),
\end{align*}
and defining ${\bf y}^{(s)}$ similarly, we note that
$$\frac {G^{(s)}({\bf x}^{(s)},{\bf y}^{(s)})} {[{\bf x}^{(s)}]^{R^{(s)}}[{\bf y}^{(s)}]^{C^{(s)}}} \quad=\quad
\prod_{k=1}^s \prod_{\ell=1}^s \frac {G({\bf x}^k, {\bf y}^\ell)} {({\bf x}^k)^R ({\bf y}^\ell)^C}.$$
From this it follows that
$$\frac{1}{s^2} \ln\left(\inf_{x^k_i,y^\ell_j>0} \frac {G^{(s)}({\bf x}^{(s)},{\bf y}^{(s)})} {[{\bf x}^{(s)}]^{R^{(s)}}[{\bf y}^{(s)}]^{C^{(s)}}} \right) \quad=\quad
\ln\left(\inf_{x_i,y_j>0} \frac{G({\bf x},{\bf y})}{{\bf x}^R{\bf y}^C} \right)$$
for all $s\geq 1$.

Inspection of the formulas for $\ln{\mathcal N}$ and $\heartsuit$ shows that both of these terms from \eqref{errorbound} have growth of order $O(s\ln s)$ as $s\rightarrow\infty$. Therefore, by \eqref{errorbound},
$$\frac{1}{s^2} \ln T_{K^{(s)}}(R^{(s)},C^{(s)}) \quad=\quad \ln\left(\inf_{x_i,y_j>0} \frac{G({\bf x},{\bf y})}{{\bf x}^R{\bf y}^C} \right) + O\left(\frac{\ln s}{s}\right),$$
from which Theorem~\ref{genfunktheorem} follows. $\blacksquare$

\section{Entropy-based estimates for $T_K(R,C)$}
\label{maxent}
In section~\ref{maxentintro}, we introduced the functions $H^{\rm max}_\kappa(x)$ and $p(x;\kappa),q(x;\kappa)$; refer to equations \eqref{implicit1}, \eqref{implicit2} for an implicit description of the latter. We now list a few useful facts about $H^{\rm max}_\kappa$:

\begin{Lem}
\label{hmaxproperties}
Let $p=p(x;\kappa)$, $q=q(x;\kappa)$.
\begin{enumerate}[{\rm (i)}]
\item $H^{\rm max}_\kappa$ is strictly concave on its domain.
\item $H^{\rm max}_\kappa(x) = -[\ln p+x\ln q]$.
\item For $0<x<\kappa$,~~$\frac{d}{dx} H^{\rm max}_\kappa(x) = -\ln q$.
\end{enumerate}
\end{Lem}

{\it Proof.}\quad
First we prove claim (i). Let $x,y\in[0,\kappa]$ and $\alpha,\beta>0$ such that $\alpha+\beta=1$. We wish to prove that
$$H^{\rm max}_\kappa(\alpha x+\beta y) > \alpha H^{\rm max}_\kappa(x) + \beta H^{\rm max}_\kappa(y).$$
Let $X$ and $Y$ be independent random variables with distributions $TG(x;\kappa)$ and $TG(y;\kappa)$, respectively. Define a random variable $Z$ whose distribution is a {\sl mixture} of $X$ and $Y$ with weights $\alpha$ and $\beta$; that is,
$${\bf Pr}[Z=t] = \alpha p(x;\kappa)q(x;\kappa)^t + \beta p(y;\kappa)q(y;\kappa)^t \quad\quad{\rm for}~t=0,1,\ldots,\kappa.$$
Then
$${\bf E}[Z] = \alpha x+\beta y$$
and
$$H[Z] > \alpha H[X] + \beta H[Y]$$
(where $H$ denotes the entropy, which is well-known to be strictly concave with respect to mixture). But
$$H^{\rm max}_\kappa(\alpha x+\beta y) \geq H[Z],$$
since $H^{\rm max}_\kappa(\alpha x+\beta y)$ is the maximum entropy achieved by any random variable supported on $\{0,1,\ldots,\kappa\}$ with expectation $\alpha x+\beta y$. This concludes the proof of (i).

Claim (ii) follows readily from equations \eqref{implicit1}, \eqref{implicit2} and the properties of logarithms. Let $p=p(x;\kappa)$, $q=q(x;\kappa)$. Then by definition of entropy, we have
\begin{align*}
H^{\rm max}_\kappa(x) &= -[p\ln p + pq\ln(pq) + pq^2\ln(pq^2) + \cdots + pq^\kappa\ln(pq^\kappa)] \\
&= -[p\ln p + pq(\ln p+\ln q) + pq^2(\ln p+2\ln q) + \cdots + pq^\kappa(\ln p+\kappa\ln q)] \\
&= -[(p+pq+pq^2+\cdots+pq^\kappa)(\ln p) + (pq+2pq^2+\cdots+\kappa pq^\kappa)(\ln q)] \\
&= -[\ln p+x\ln q].
\end{align*}

Differentiating this formula with respect to $x$, and again using equations \eqref{implicit1} and \eqref{implicit2}, we obtain
\begin{align*}
(H^{\rm max}_\kappa)'(x) &= -\frac{p'}{p} - x\cdot\frac{q'}{q} -\ln q \\
&= p\cdot\left(\frac{1}{p}\right)' - p(q+2q^2+\cdots+\kappa q^\kappa)\cdot\frac{q'}{q} -\ln q \\
&= p\cdot\left(\frac{1}{p}\right)' - pq'(1+2q+\cdots+\kappa q^{\kappa-1}) -\ln q \\
&= p\cdot\left(\frac{1}{p}\right)' - p\cdot\left(\frac{1}{p}\right)' -\ln q \\
&= -\ln q.
\end{align*}
This proves claim (iii). $\square$

Like the ``$(\kappa+1)$-nomial coefficients'' of section~\ref{indepheu}, the functions $H^{\rm max}_\kappa$, $p$, and $q$ admit simple explicit formulas only when $\kappa=1$ or $\kappa=\infty$. To wit:

\begin{align*}
& H^{\rm max}_1(x)=-x\ln x-(1-x)\ln(1-x) &&~~\quad p(x;1)=1-x &&~~\quad q(x;1)=\frac{x}{1-x} \\
& H^{\rm max}_\infty(x)=(x+1)\ln(x+1)-x\ln x &&~~\quad p(x;\infty)=\frac{1}{x+1} &&~~\quad q(x;\infty)=\frac{x}{x+1}
\end{align*}

In fact, there is a close relationship between these functions and the $(\kappa+1)$-nomial coefficients:
\begin{Lem}
\label{knomial}
Let $\kappa\in{\mathbb Z}_{>0}\cup\{\infty\}$. Let $n,r$ be integers ($n>0,~0\leq r\leq n\kappa$). Then
$$\lim_{s\rightarrow\infty} \frac{1}{s}\ln{sn\choose sr}_\kappa \quad=\quad nH^{\rm max}_\kappa\left(\frac{r}{n}\right).$$
\end{Lem}

{\it Proof.}\quad
Let $X_1,X_2,\ldots$ be independent random variables, each with distribution $TG\left(\frac{r}{n};\kappa\right)$. Let $X=(X_1,\ldots,X_{sn})$.

Observe that if ${\bf x},{\bf x}'\in\{0,1,\ldots,\kappa\}^s$, then
$$\frac{{\bf Pr}[X={\bf x}']}{{\bf Pr}[X={\bf x}]} \quad=\quad q^{|{\bf x}'|-|{\bf x}|}$$
(where $|{\bf x}| := \sum_{i=1}^s x_i$). In particular, all values of $X$ with equal sum of coordinates are equiprobable. Let ${\bf x}_*$ denote an arbitrary value for $X$ satisfying $|{\bf x}_*|=sr$.

Recall that the {\sl Shannon self-information} of a value $X={\bf x}$ is defined as
$$I({\bf x}) := -\ln {\bf Pr}[X={\bf x}];$$
the entropy of $X$ is the expected self-information of its value. Thus we have
\begin{align}
snH^{\rm max}_\kappa\left(\frac{r}{n}\right) \quad=\quad H[X] \quad&=\quad {\bf E}[I(X)] \nonumber\\
&=\quad I({\bf x}_*) - (\ln q)~{\bf E}\big[|X|-sr\big] \nonumber\\
&=\quad I({\bf x}_*) \nonumber\\
&=\quad -\ln\left[{sn\choose sr}_\kappa^{-1}\cdot {\bf Pr}\big[|X|=sr\big]\right] \nonumber\\
&=\quad \ln{sn\choose sr}_\kappa - \ln {\bf Pr}\big[|X|=sr\big]. \label{dependencebound}
\end{align}

Note that the probability mass function for each $X_i$ is log-concave on ${\mathbb Z}$. We apply a local limit theorem of Bender (see Appendix, Theorem 7) using
\begin{align*}
&\zeta_p=X_1+\cdots+X_p,
&&\sigma_p^2=p\cdot\var(X_1),
&&\mu_p=p\cdot\frac{r}{n},
&&{\rm and} &&x=0,
\end{align*}
with the normality hypothesis secured via Lyapunov's central limit theorem (Appendix, Theorem 6), to infer
$$\lim_{p\rightarrow\infty} \sigma_p {\bf Pr}\left[\zeta_p = \left\lfloor p\cdot\frac{r}{n}\right\rfloor\right] \quad=\quad \frac{1}{\sqrt{2\pi}},$$
and thence
$${\bf Pr}\big[|X|=sr\big] \quad\sim\quad (2\pi sn\var(X_1))^{-1/2} \quad=\quad \Theta(s^{-1/2}).$$
Substituting into \eqref{dependencebound}, we conclude that
$$snH^{\rm max}_\kappa\left(\frac{r}{n}\right) \quad=\quad \ln{sn\choose sr}_\kappa - \Theta(\ln s);$$
in particular, this proves the lemma. $\square$

\subsection{A dual to the optimization problem in Theorem~\ref{genfunktheorem}.}
In Theorem~\ref{genfunktheorem}, we ``computed'' $\ln T_K(R,C)$ asymptotically as
\begin{equation}
\ln\left(\inf_{x_i,y_j>0} \frac{G({\bf x},{\bf y})}{{\bf x}^R{\bf y}^C}\right). \label{genfunklimit}
\end{equation}
We now show:
\begin{Lem}
\label{duality}
Suppose $k_{ij} = \kappa$ for $1\leq i\leq m$, $1\leq j\leq n$. Then the quantity given in \eqref{genfunklimit} is equal to
\begin{equation}
\max_{Z\in \Pi_\kappa(R,C)} \sum_{i=1}^m \sum_{j=1}^n H^{\rm max}_\kappa(z_{ij}). \label{entropylimit}
\end{equation}
\end{Lem}

{\it Proof.}\quad
By Lemma~\ref{hmaxproperties}, $H^{\rm max}_\kappa(x)$ is strictly concave. Also, $(H^{\rm max}_\kappa)'(x) = \ln q(x;\kappa)$ approaches $\infty$ as $x\rightarrow 0$ and $-\infty$ as $x\rightarrow\kappa$; therefore, the maximum in \eqref{entropylimit} is well-defined and is attained in the relative interior of $\Pi_\kappa(R,C)$. For the remainder of this proof, let $Z$ denote the (unique) location at which the maximum is attained, and let $p_{ij}:=p(z_{ij};\kappa)$, $q_{ij}:=q(z_{ij};\kappa)$.

Since $Z$ is in the interior of $\Pi_\kappa(R,C)$, the local defining equations for $\Pi_\kappa(R,C)$ at $Z$ are just
$$\sum_{j=1}^n a_{ij} = r_i~~(1\leq i\leq m)\quad{\rm and}\quad\sum_{i=1}^m a_{ij} = c_j~~(1\leq j\leq n).$$
Introducing Lagrange multipliers for these constraints, we infer that $\ln q_{ij} = \lambda_i+\mu_j$ for some constants $\lambda_1,\ldots,\lambda_m,\mu_1,\ldots,\mu_n$. Define $\xi_i:=e^{\lambda_i}$, $\eta_j=e^{\mu_j}$; thus $q_{ij}=\xi_i\eta_j$. Dividing equation \eqref{implicit2} by equation \eqref{implicit1} (see section~\ref{maxentintro}), we obtain
$$z_{ij} = \frac{\xi_i\eta_j+2(\xi_i\eta_j)^2+\cdots+\kappa(\xi_i\eta_j)^\kappa} {1+\xi_i\eta_j+(\xi_i\eta_j)^2+\cdots+(\xi_i\eta_j)^\kappa}.$$

For real-valued ${\bf t}=(t_1,\ldots,t_m)$ and ${\bf s}=(s_1,\ldots,s_n)$, let
\begin{align*}
\psi({\bf t},{\bf s}) ~~&:=~~ \ln \frac{G({\bf x},{\bf y})}{{\bf x}^R{\bf y}^C} \bigg|{}_{\substack{x_i=e^{t_i} \\ y_j=e^{s_j}}} \\
&=~~ -\sum_{i=1}^m r_it_i -\sum_{j=1}^n c_js_j + \sum_{i=1}^m\sum_{j=1}^n \ln\left(1+e^{t_i+s_j}+e^{2(t_i+s_j)}+\cdots+e^{\kappa(t_i+s_j)}\right).
\end{align*}
This function is strictly convex, and has a critical point (hence a global minimum) at $({\bf t},{\bf s})$ if and only if the gradient is zero, that is, if
\begin{align*}
r_i ~~&=~~ \sum_{j=1}^n \frac{e^{t_i+s_j}+2e^{2(t_i+s_j)}+\cdots+\kappa e^{\kappa(t_i+s_j)}} {1+e^{t_i+s_j}+e^{2(t_i+s_j)}+\cdots+e^{\kappa(t_i+s_j)}}, &&1\leq i\leq m \\
{\rm and}~~c_j ~~&=~~ \sum_{i=1}^m \frac{e^{t_i+s_j}+2e^{2(t_i+s_j)}+\cdots+\kappa e^{\kappa(t_i+s_j)}} {1+e^{t_i+s_j}+e^{2(t_i+s_j)}+\cdots+e^{\kappa(t_i+s_j)}}, &&1\leq j\leq n.
\end{align*}
These conditions are satisfied at ${\bf t}=(\lambda_1,\ldots,\lambda_m)$ and ${\bf s}=(\mu_1,\ldots,\mu_n)$. The minimum value of $\psi$ is thus
\begin{align*}
\psi({\bf t},{\bf s}) ~~&=~~
-\sum_{i=1}^m r_i\lambda_i -\sum_{j=1}^n c_j\mu_j + \sum_{i=1}^m \sum_{j=1}^n \ln\left(1+\xi_i\eta_j+(\xi_i\eta_j)^2+\cdots+(\xi_i\eta_j)^\kappa\right) \\
&=~~ \sum_{i=1}^m \sum_{j=1}^n \left[-z_{ij}(\lambda_i+\mu_j) + \ln(1+q_{ij}+q_{ij}^2+\cdots+q_{ij}^\kappa)\right] \\
&=~~ \sum_{i=1}^m \sum_{j=1}^n \left[-z_{ij}\ln q_{ij} + \ln\left(\frac{1}{p_{ij}}\right)\right] \\
&=~~ \sum_{i=1}^m \sum_{j=1}^n H^{\rm max}_\kappa(z_{ij}).
\end{align*}
This proves the lemma. $\square$

\subsection{Proof of Theorem~\ref{entropytheorem}}
\label{hproof}
Part (i) follows directly from Theorem~\ref{genfunktheorem} and Lemma~\ref{duality}.

For part (ii), recall that
$$I_\kappa(R,C) = {mn\choose N}_\kappa^{-1} ~\prod_{i=1}^m {n\choose r_i}_\kappa~\prod_{j=1}^n {m\choose c_j}_\kappa$$
(Definition~\ref{indepest}). Thus
$$I_\kappa(R^{(s)},C^{(s)}) = {s^2mn\choose s^2N}_\kappa^{-1} ~\left(\prod_{i=1}^m {sn\choose sr_i}_\kappa\right)^s ~\left(\prod_{j=1}^n {sm\choose sc_j}_\kappa\right)^s.$$
Applying Lemma~\ref{knomial}, we obtain
\begin{align*}
\ln I_\kappa(R^{(s)},&C^{(s)}) ~~=~~ -\left[s^2mn H^{\rm max}_\kappa\left(\frac{N}{mn}\right) + o(s^2)\right]
+ s\sum_{i=1}^m \left[snH^{\rm max}_\kappa\left(\frac{r_i}{n}\right)+o(s)\right] \\
&\hskip 100 pt + s\sum_{j=1}^n \left[smH^{\rm max}_\kappa\left(\frac{c_j}{m}\right)+o(s)\right] \\
&=~~ s^2\left[-mnH^{\rm max}_\kappa \left(\frac{N}{mn}\right) + n\sum_{i=1}^m H^{\rm max}_\kappa \left(\frac{r_i}{n}\right)
+ m\sum_{j=1}^n H^{\rm max}_\kappa \left(\frac{c_j}{m}\right) + o(1)\right].
\end{align*}

This completes the proof of Theorem~\ref{entropytheorem}. $\blacksquare$

\subsection{The entropy loss function}
\label{entropyloss}
The following function plays a key role in the proof of Theorem~\ref{attraction}:
\begin{Def}
\label{hloss}
Fix $\kappa\in{\mathbb Z}_{\geq 0}\cup\{\infty\}$.
Given nonnegative $\alpha_1,\alpha_2,\ldots,\alpha_n$ such that $\alpha_1+\alpha_2+\cdots+\alpha_n=1$, let
$$J(r) ~~=~~ J_{\alpha,\kappa}(r) ~~:=~~
nH^{\rm max}_\kappa\left(\frac{r}{n}\right) - \sum_{j=1}^n H^{\rm max}_\kappa(r\alpha_j)$$
for all $r\geq 0$ such that $r\alpha_1,r\alpha_2,\ldots,r\alpha_n\leq\kappa$.
\end{Def}

In the spirit of the remarks at the beginning of section~\ref{maxentintro}, we propose the following (informal) interpretation for $J(r)$. (The reader who is interested only in formal proof may skip to section~\ref{attractionproof}.)

Suppose $A$ is an $m\times n$ contingency table, about which we know only that the sum of entries is $N$. In order to guess what $A$ looks like, we might sample $mn$ entries independently from $TG\left(\frac{N}{mn};\kappa\right)$. The resulting matrix might not have sum of entries exactly equal to $N$, but at least it is correct in expectation, and all actual tables with sum of entries $N$ are equally likely to be chosen. The entropy of our random model is $mnH^{\rm max}_\kappa\left(\frac{N}{mn}\right)$.

If we subsequently learn that $A$ has column sums $C=(c_1,\ldots,c_n)$, then we might revise our model by drawing entries in the $j^{\rm th}$ column from $TG\left(\frac{c_j}{m};\kappa\right)$. The entropy of the random model is then $m\sum_{j=1}^n H^{\rm max}_\kappa\left(\frac{c_j}{m}\right)$. The information gained may be measured by the entropy lost, which is equal to $mJ\left(\frac{N}{m}\right)$ for $\alpha_j = \frac{c_j}{N}$\break
$(j=1,\ldots,n)$. 

The same comparison may be made in the presence of known row sums $R=(r_1,\ldots,r_m)$. Before learning $C$, we sample entry $(i,j)$ from $TG\left(\frac{r_i}{n};\kappa\right)$; after learning $C$, we sample entry $(i,j)$ from $TG\left(\frac{r_ic_j}{N};\kappa\right)$. Of course, this guess is entirely na\"ive---it essentially regards the rank~1 table as typical, which is a doubtful assumption. We do not propose that this is the {\it best} guess, but it does have the advantage of being computable. The entropy lost to this model when incorporating $C$ in the presence of known $R$ is
$$\sum_{i=1}^m J(r_i),$$
where $\alpha_j = \frac{c_j}{N}$~~$(j=1,\ldots,n)$. Note that this quantity is defined if and only if the rank~1 table has all entries~$\leq\kappa$.

We would like to know if the row margin $R$ and the column margin $C$ are positively correlated. Intuitively, this means that the revelation of $C$ produces less {\sl surprise} (entropy loss) when $R$ is known in advance than when $R$ is not known. In our random model, this is expressed by the inequality
\begin{equation}
mJ\left(\frac{N}{m}\right) \quad\geq\quad J(r_1)+J(r_2)+\cdots+J(r_m). \label{hlossinequality}
\end{equation}
Now let us see if we can make a formal argument out of the intuition.

Theorem~\ref{entropytheorem}, part (i), implies that
$$\lim_{s\rightarrow\infty} \frac{1}{s^2} \ln T_\kappa(R^{(s)},C^{(s)}) \quad\geq\quad
\sum_{i=1}^m \sum_{j=1}^n H^{\rm max}_\kappa\left(\frac{r_ic_j}{N}\right);$$
again, the rank~1 matrix is chosen here purely out of convenience (we could have substituted any $Z\in\Pi_\kappa(R,C)$).
Together with part (ii), this means that in order to prove
$$\lim_{s\rightarrow\infty} \frac{1}{s^2} \ln T_\kappa(R^{(s)},C^{(s)}) \quad\geq\quad
\lim_{s\rightarrow\infty} \frac{1}{s^2} \ln I_\kappa(R^{(s)},C^{(s)}),$$

it suffices to show that
$$mnH^{\rm max}_\kappa\left(\frac{N}{mn}\right)
\quad\geq\quad
- \sum_{i=1}^m \sum_{j=1}^n H^{\rm max}_\kappa\left(\frac{r_ic_j}{N}\right)
+ m\sum_{j=1}^n H^{\rm max}_\kappa\left(\frac{c_j}{m}\right)
+ n\sum_{i=1}^m H^{\rm max}_\kappa\left(\frac{r_i}{n}\right)$$
or, equivalently (by rearranging terms),
$$mnH^{\rm max}_\kappa\left(\frac{N}{mn}\right)
- m\sum_{j=1}^n H^{\rm max}_\kappa\left(\frac{c_j}{m}\right)
\quad\geq\quad
n\sum_{i=1}^m H^{\rm max}_\kappa\left(\frac{r_i}{n}\right)
- \sum_{i=1}^m \sum_{j=1}^n H^{\rm max}_\kappa\left(\frac{r_ic_j}{N}\right).$$

This is precisely what is asserted by the inequality \eqref{hlossinequality}. In order to prove Theorem~\ref{attraction}, we will show that under its hypotheses, \eqref{hlossinequality} is true (and meaningful, i.e.~$A^{\rm ind}$ has entries~$\leq\kappa$).

\subsection{Proof of Theorem~\ref{attraction}}
\label{attractionproof}
Assume the hypotheses of Theorem~\ref{attraction}; note in particular that $\kappa\geq 2$. Let
$$\alpha_j ~~:=~~ \frac{c_j}{N} \quad\quad (1\leq j\leq n).$$

Consider the function
$$\phi(x) \quad:=\quad x^2(H^{\rm max}_\kappa)''(x) \quad=\quad -x^2\cdot\frac{q'(x;\kappa)}{q(x;\kappa)}$$
(all derivatives being with respect to $x$; the second equality follows from Lemma~\ref{hmaxproperties}(iii)).

The above formula defines $\phi(x)$ only for $0\leq x\leq\kappa$, but we claim that $\phi(x)$ can be extended analytically to a neighborhood of $x=0$.

{\it Proof of claim:}\quad
Equations \eqref{implicit1} and \eqref{implicit2} (section~\ref{maxentintro}) yield
$$x = \frac{q+2q^2+\cdots+\kappa q^\kappa}{1+q+q^2+\cdots+q^\kappa},$$
where $q=q(x;\kappa)$. Although this formula has only been assigned meaning for $q\geq 0$, it shows that $x$ (as a function of $q$) can be extended analytically to a neighborhood of $q=0$; the Maclaurin series is $x=q+q^2+O(q^3)$. Since $\frac{dx}{dq}\neq 0$ at $q=0$, it follows that the inverse function $q(x;\kappa)$ is also defined and analytic in a neighborhood of $x=0$, with Maclaurin series $q=x-x^2+O(x^3)$. Applying l'H\^opital's rule, we see that the singularity of $\phi$ at $x=0$ is removable, so $\phi(x)$ is locally analytic there. $\square$

We compute the Maclaurin series of $\phi(x)$:
$$\phi(x) \quad=\quad -x\cdot\frac{1-2x+O(x^2)}{1-x+O(x^2)} \quad=\quad -x+x^2+O(x^3).$$
Since the coefficient of $x^2$ is positive, $\phi(x)$ is strictly convex in a neighborhood of $x=0$. Choose $\delta\in(0,1)$ such that $\phi(x)$ is strictly convex in the interval\break
$|x|\leq\delta\kappa$.

Because $\delta<1$,~~$J(r)$ (see definition~\ref{hloss}) is defined and differentiable at $r=r_1,\ldots,r_m$. Differentiating, we have
$$J'(r) ~~=~~ (H^{\rm max}_\kappa)'\left(\frac{r}{n}\right) - \sum_{j=1}^n \alpha_j(H^{\rm max}_\kappa)'(r\alpha_j)$$
and
\begin{align*}
J''(r) ~~&=~~ \frac{1}{n}(H^{\rm max}_\kappa)''\left(\frac{r}{n}\right) - \sum_{j=1}^n \alpha_j^2(H^{\rm max}_\kappa)''(r\alpha_j) \\
&=~~ \frac{n}{r^2}\phi\left(\frac{r}{n}\right) - \sum_{j=1}^n \frac{1}{r^2}\phi(r\alpha_j).
\end{align*}
By the convexity of $\phi(x)$, we have $J''(r)\leq 0$ for $0<r\leq \frac{\delta\kappa}{\max\{\alpha_1,\ldots,\alpha_n\}}$; the inequality is strict if $\alpha_1,\ldots,\alpha_n$ are not all equal. Therefore, $J(r)$ is concave on (the closure of) that interval, and strictly concave if $\alpha_1,\ldots,\alpha_n$ are not all equal. By our assumption that $r_ic_j\leq\delta\kappa N$, it follows that $r_1,\ldots,r_m$ are in that interval.

Thus, inequality \eqref{hlossinequality} holds:
$$mJ\left(\frac{N}{m}\right) \quad\geq\quad J(r_1)+J(r_2)+\cdots+J(r_m),$$
with strict inequality if $\alpha_1,\ldots,\alpha_n$ are not all equal and $r_1,\ldots,r_m$ are also not all equal.

When the function $J$ is evaluated throughout this inequality, we obtain
$$mnH^{\rm max}_\kappa\left(\frac{N}{mn}\right)
- m\sum_{j=1}^n H^{\rm max}_\kappa\left(\frac{c_j}{m}\right)
\quad\geq\quad
n\sum_{i=1}^m H^{\rm max}_\kappa\left(\frac{r_i}{n}\right)
- \sum_{i=1}^m \sum_{j=1}^n H^{\rm max}_\kappa\left(\frac{r_ic_j}{N}\right).$$
As explained at the end of section~\ref{entropyloss}, this implies the last link, and Theorem~\ref{entropytheorem} implies the first, in this chain of inequalities:
\begin{align*}
\lim_{s\rightarrow\infty} \frac{1}{s^2} \ln T_\kappa(R^{(s)},C^{(s)}) \quad&\geq\quad
\max_{Z\in\Pi_\kappa(R,C)} \sum_{i=1}^m \sum_{j=1}^n H^{\rm max}_\kappa(z_{ij}) \\
&\geq\quad
\sum_{i=1}^m \sum_{j=1}^n H^{\rm max}_\kappa\left(\frac{r_ic_j}{N}\right) \\
&\geq\quad
\lim_{s\rightarrow\infty} \frac{1}{s^2} \ln I_\kappa(R^{(s)},C^{(s)}).
\end{align*}

If $\alpha_1,\ldots,\alpha_n$ are not all equal and $r_1,\ldots,r_m$ are not all equal, then the last inequality in this chain is strict. This completes the proof of Theorem~\ref{attraction}. $\blacksquare$

\subsection{Prospects for negative correlation of margins}
\label{repulsion}
Recall that for $\kappa=1$, all pairs of margins $(R,C)$ have either zero or negative asymptotic correlation (specifically, negative unless either $R$ or $C$ is a constant vector). For $\kappa=\infty$, the sign of correlation is reversed. We expect that these are the only ``pure'' cases: that is, when $1<\kappa<\infty$, there are some positively correlated pairs of margins as well as some negatively correlated pairs. Theorem~\ref{attraction} asserts half of this conjecture: for $\kappa\geq 2$, any sufficiently sparse margins are asymptotically positively correlated. (By symmetry, ``co-sparse margins''---those which force most entries to be close to $\kappa$---are also positively correlated.)

Numerical evidence and heuristic arguments suggest that, for all $\kappa<\infty$, margins which are neither sparse nor co-sparse---or more specifically, close to $R=(\frac{n\kappa}{2},\ldots,\frac{n\kappa}{2})$ and $C=(\frac{m\kappa}{2},\ldots,\frac{m\kappa}{2})$---are negatively correlated. For example, we have used Theorem~\ref{entropytheorem} to compute
\begin{equation}
\lim_{s\rightarrow\infty} \frac{1}{s^2} \left[\ln T_\kappa(R^{(s)},C^{(s)}) - \ln I_\kappa(R^{(s)},C^{(s)})\right] \label{heuristicerror}
\end{equation}
for margins of the form $R=C=(\gamma,\gamma+\varepsilon,\gamma+2\varepsilon,\ldots,\gamma+(n-1)\varepsilon)$, $n=2,3,4,5$. In these tests, when $\kappa=2$ and $\varepsilon=.02$, expression \eqref{heuristicerror} turns out to be negative (indicating negative correlation of the margins) roughly when $.09n<\gamma<1.89n$, suggesting a threshold of $\delta\approx .05$. When $\kappa=10$ and $\varepsilon=.1$, expression \eqref{heuristicerror} is negative roughly when $1.5n<\gamma<8.5n$, suggesting $\delta\approx .15$. As $\kappa\rightarrow\infty$, the sharp value of $\delta$ in Theorem~\ref{attraction} appears to grow, but a threshold remains.

An intuitive gloss on this phenomenon is that the distribution $TG(x;\kappa)$ ``looks like'' a geometric distribution when $x\approx 0$ (or $x\approx\kappa$), but looks more like a Bernoulli distribution when $x$ is at neither extreme. In the former case, the ``lid'' $\kappa$ (or the floor 0) is remote from typical values, so the behavior observed when $\kappa=\infty$ dominates. In the latter case, the $\kappa=1$ behavior seems to dominate.

The fundamental difference between these cases is suggested by the function $\phi(x)$ in the proof of Theorem~\ref{attraction}. When $\kappa=\infty$, this function is convex throughout its domain; when $\kappa=1$, it is concave; and when $1<\kappa<\infty$, this function is convex near the origin, but has an inflection point.

We can show that $\phi(x)$ is concave for $x\approx\frac{\kappa}{2}$, so what are the obstacles to a reverse Theorem~\ref{attraction}? There are two. In the proof of Theorem~\ref{attraction}, we relied on the fact that
\begin{equation}
\max_{Z\in\Pi_\kappa(R,C)} \sum_{i=1}^m \sum_{j=1}^n H^{\rm max}_\kappa(z_{ij})
\quad\geq\quad
\sum_{i=1}^m \sum_{j=1}^n H^{\rm max}_\kappa\left(\frac{r_ic_j}{N}\right), \label{maxinequality}
\end{equation}
allowing us to use the rank 1 matrix $A^{\rm ind}=\left(\frac{r_ic_j}{N}\right)$ as a proxy for the unknown $Z$ which achieves the maximum. This matrix does not necessarily have entries $\leq\kappa$; we were able to assume that it does only because our assumption of sparse margins did double duty. This is the first obstacle to a reverse Theorem~\ref{attraction}; the second is that, even if we could find another plausible candidate for $Z$, we could not make an assumption like \eqref{maxinequality} with the inequality reversed. Thus, to prove a negative correlation result, we believe it is crucial to understand something about {\it where} the maximum on the left-hand side of \eqref{maxinequality} is achieved.

\section{Appendix}
This section contains some matter which was deferred from earlier sections.

\subsection{Proof of Lemma~\ref{fproperties}}
Claim (i) is trivial, since we can set $p=1$, $\alpha_1=1$, $R^1=R$, $C^1=C$ in Definition~\ref{f}.

For claim (ii), it suffices to show that if $\alpha+\beta=1$, then
\begin{equation}
\alpha f(R^1,C^1) + \beta f(R^2,C^2) \leq f(\alpha R^1+\beta R^2,\alpha C^1+\beta C^2). \label{concavity}
\end{equation}
By Definition~\ref{f}, there exist $\gamma_1,\ldots,\gamma_p\geq 0$; $R^{11},\ldots,R^{1p}$; and $C^{11},\ldots,C^{1p}$ such that
\begin{align*}
&\sum_{t=1}^p \gamma_t=1, \quad\quad
\sum_{t=1}^p \gamma_tR^{1t} = R, \quad\quad
\sum_{t=1}^p \gamma_tC^{1t} = C, \\
&{\rm and}\quad
f(R^1,C^1) = \sum_{t=1}^p \gamma_t\ln T_K(R^{1t},C^{1t}).
\end{align*}

Likewise, there exist $\delta_1,\ldots,\delta_q\geq 0$; $R^{21},\ldots,R^{2q}$; and $C^{21},\ldots,C^{2q}$ such that
\begin{align*}
&\sum_{t=1}^q \delta_t=1, \quad\quad
\sum_{t=1}^q \delta_tR^{2t} = R, \quad\quad
\sum_{t=1}^q \delta_tC^{2t} = C, \\
&{\rm and}\quad
f(R^2,C^2) = \sum_{t=1}^q \delta_t\ln T_K(R^{2t},C^{2t}).
\end{align*}

Note that
\begin{align*}
&\sum_{t=1}^p \alpha\gamma_t + \sum_{t=1}^q \beta\delta_t = 1, \quad\quad
&\sum_{t=1}^p \alpha\gamma_tR^{1t} + \sum_{t=1}^q \beta\delta_tR^{2t} = \alpha R^1+\beta R^2, \\
&~ &{\rm and}\quad\quad
\sum_{t=1}^p \alpha\gamma_tC^{1t} + \sum_{t=1}^q \beta\delta_tC^{2t} = \alpha C^1+\beta C^2;
\end{align*}
applying Definition~\ref{f} to $f(\alpha R^1+\beta R^2,\alpha C^1+\beta C^2)$, we obtain equation \eqref{concavity} and thus claim (ii).

It is clear that $f$ is defined only on the convex hull of all $(R,C)$ for which $T_K(R,C)>0$; this region is a subset of $\Pi_K(R,C)$, proving claim (iii). $\square$

\subsection{Limit theorems of Lyapunov and Bender}
We make use of the following theorems in the proof of Lemma~\ref{knomial}:

\begin{Thm}[Lyapunov's central limit theorem]
Suppose that $\big(X_n: n\in{\mathbb N}\big)$ is a sequence of independent random variables, such that $\mu_n:={\bf E}[X_n]$ and $\sigma_n^2:=\var[X_n]$ are finite.
Let $\zeta_n=X_1+\cdots+X_n$, and define $m_n:={\bf E}[\zeta_n]=\mu_1+\cdots+\mu_n$, $s_n^2:=\var[\zeta_n]=\sigma_1^2+\cdots+\sigma_n^2$. If
$$\lim_{n\rightarrow\infty} \frac{1}{s_n^{2+\delta}} \sum_{k=1}^n {\bf E}\left[|X_k-\mu_k|^{2+\delta}\right] = 0$$
for some $\delta>0$, then
$$\lim_{n\rightarrow\infty} {\bf Pr}\left[ \zeta_n < s_nx+m_n \right] = \frac{1}{\sqrt{2\pi}} \int_{-\infty}^x e^{-t^2/2} dt$$
for all $x\in{\mathbb R}$.
\end{Thm}

\begin{Thm}[Bender local limit theorem]
Suppose that $\big(\zeta_n: n\in{\mathbb N}\big)$ is a sequence of integer-valued random variables and $\big(\sigma_n\big)$ and $\big(\mu_n\big)$ are sequences of real numbers, such that
$$\lim_{n\rightarrow\infty} {\bf Pr}\left[ \zeta_n < \sigma_nx+\mu_n \right] = \frac{1}{\sqrt{2\pi}} \int_{-\infty}^x e^{-t^2/2} dt$$
for all $x\in{\mathbb R}$. Also suppose that $\sigma_n\rightarrow\infty$ as $n\rightarrow\infty$. Further, suppose that, for every $n$, the sequence $b_n(t):={\bf Pr}(\zeta_n=t)$ is properly log-concave with respect to $t$. Then
$$\lim_{n\rightarrow\infty} \sigma_n{\bf Pr}\left[\zeta_n = \lfloor\sigma_n x+\mu_n\rfloor\right] = \frac{1}{\sqrt{2\pi}} e^{-x^2/2}$$
uniformly for all $x\in{\mathbb R}$.
\end{Thm}

The Lyapunov theorem is well-known; for a proof, see e.g.~\cite{Skorokhod2005}. The Bender theorem originally appeared in~\cite{Bender1973}; see also~\cite{Engel1997}, which states and proves the theorem in a form more similar to the above.

\section{Acknowledgments}
The author thanks Alexander Barvinok for posing the problems addressed by this paper, and for many useful discussions.

\bibliographystyle{plain}
\bibliography{contingency_bibliography}

\end{document}